\documentclass[14pt,twoside]{extarticle}

\pdfoutput=1 

\usepackage{cmap}
\usepackage[utf8]{inputenc}
\usepackage[english]{babel}

\usepackage{amsmath,amsfonts,amsthm,amssymb}
\usepackage[margin=1in,marginparsep=22pt]{geometry}

\usepackage[unicode]{hyperref}
\usepackage{caption}
\usepackage{float}
\usepackage{wrapfig}
\usepackage{enumitem}
\usepackage{changepage}
\usepackage{mathtools}
\usepackage{etoolbox}
\usepackage{graphicx}
\usepackage{xstring}
\usepackage{xspace}

\newbool{online}
\newbool{maketoc} 

\booltrue{online}\booltrue{maketoc}

\ifonline
\else
	\hypersetup{nolinks=true}%
\fi

\newcommand{\tsc}{\scshape\small}

\newcommand{\txtsc}[1]{{\textsc{#1}}}

\newcommand\headtitle{Headtitle}
\newcommand\subtitle{Subtitle}

\newcommand\sectitle{\headtitle}
\newcommand\mauthor{Author}

\newcommand{\vvs}{\vspace{5pt}}
\newcommand{\vvm}{\vspace{10pt}}

\newlength\gvhs
\newlength\gvls
\newlength\gvms
\newlength\gvss
\newlength\gvts
\newlength\gvxs
\newcommand{\vmskip}{\vspace{\gvms}}
\newcommand{\vsskip}{\vspace{\gvss}}
\newcommand{\vtskip}{\vspace{\gvts}}

\newenvironment{mcenter}%
{\parskip=0pt\par\nopagebreak\centering}%
{\par\noindent\ignorespacesafterend}%

\newenvironment{qte}[1][-2ex]{\vsskip\pagebreak[3]\begingroup\wpic{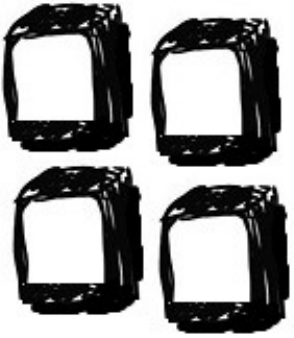}{w=0.06,s=0.15,p=l,l=99,v=#1}\small}%
{\par\endgroup\vsskip}

\newcommand\glm[1]{\guillemotleft\kern1pt#1\kern1pt\guillemotright}
\newcommand\tquote[1]{\noindent{\it\glm{#1}}}
\newcommand{\tbrk}[1]{$[$\kern1pt#1\kern1pt$]$}
\newcommand{\tnumsign}[1]{\textnumero\,$#1$}
\newcommand{\xxx}{\vvm\centerline{*\quad *\quad *}\vvs}

\newcommand{\pp}{\,.}

\newcommand\ifempty[1]{\def\temp{#1}\ifx\temp\empty}

\newcommand\BeginPic{\begingroup}
\newcommand\EndPic[1][0.1pt]{\vspace{#1}\endgroup\WFclear}

\newcommand{\ccaption}[1]{\begin{mcenter}{\tsc #1}\end{mcenter}\vtskip}

\newcommand{\pic}[2]{\includegraphics[scale=#1]{#2}}

\newcommand{\ccpic}[3]{\begin{figure}[H]\begin{mcenter}%
\pic{#1}{#2}\ifempty{#3}\empty\else\par{\tsc #3}\vspace{-5pt}\fi%
\end{mcenter}\end{figure}}

\makeatletter

\newcounter{wp@lines}

\define@key{wpkeys}{width}{\def\wp@width{#1}\relax}
\define@key{wpkeys}{w}{\def\wp@width{#1}\relax}
\define@key{wpkeys}{scale}{\def\wp@scale{#1}\relax}
\define@key{wpkeys}{s}{\def\wp@scale{#1}\relax}
\define@key{wpkeys}{lines}{\setcounter{wp@lines}{#1}\relax}
\define@key{wpkeys}{l}{\setcounter{wp@lines}{#1}\relax}
\define@key{wpkeys}{vskip}{\def\wp@vskip{#1}\relax}
\define@key{wpkeys}{v}{\def\wp@vskip{#1}\relax}
\define@key{wpkeys}{pos}{\def\wp@side{#1}\relax}
\define@key{wpkeys}{p}{\def\wp@side{#1}\relax}
\define@key{wpkeys}{caption}{\def\wp@caption{#1}\relax}
\define@key{wpkeys}{c}{\def\wp@caption{#1}\relax}

\newcommand{\wpic}[2]%
{%
  \def\wp@width{1.0}%
  \def\wp@scale{1.0}%
  \setcounter{wp@lines}{-1}%
  \def\wp@vskip{-2ex}%
  \def\wp@side{o}%
  \let\wp@caption\empty%
  \setkeys{wpkeys}{#2}%
  \ifx\empty\wp@caption\relax\def\cptext{\empty}\else
                             \def\cptext{\par{\txtsc{\wp@caption}}\vspace{-5pt}}\fi
  \ifnum\thewp@lines<0\relax
      \begin{wrapfigure}{\wp@side}{\wp@width\textwidth}%
      \centering{\vspace{\wp@vskip}\pic{\wp@scale}{#1}}%
      \cptext
      \end{wrapfigure}%
  \else
      \begin{wrapfigure}[\thewp@lines]{\wp@side}{\wp@width\textwidth}%
      \centering{\vspace{\wp@vskip}\pic{\wp@scale}{#1}}%
      \cptext
      \end{wrapfigure}%
  \fi
}

\makeatother

\newcounter{fsection}

\newcommand{\nsfrac}[2]{\,\raise3pt\hbox{$\scriptstyle{#1}$\hskip-2pt}/\raise-3pt\hbox{\hskip-2pt$\scriptstyle{#2}$}}

\newcommand\ftnnote[1]{\footnote{\ #1}}

\def\firstchar#1{\StrChar{#1}{1}}
\def\Fio(#1 #2.#3){#1 #2.{\kern 1pt}#3}
\def\fio(#1 #2.#3){\firstchar{#1}.\kern 1pt #2.\kern 2pt #3}
\def\fnln(#1 #2){\firstchar{#1}.\kern 2pt #2}
\def\lfp(#1 #2.#3.){#1~\firstchar{#2}.\kern 2pt\firstchar{#3}.}
\def\fior(#1.#2.#3){#1.{\kern 1pt}#2.{\kern 2pt}#3}
\def\fim(#1.#2){#1.{\kern 1pt}#2}

\newcommand\infolink[2]{\href{#2}{#1}}

\ifonline%
\else%

\renewcommand{\infolink}[2]{\textbf{#1}}
\fi

\newcommand\stpb{St.\kern0.1em Petersburg\xspace}

\renewcommand\headtitle{Math Matters of the Past}
\renewcommand\subtitle{The Very First Mathematical Olympiads}
\renewcommand\mauthor{Dmitri Fomin}
\renewcommand\sectitle{\subtitle}

\makeatletter
\def\@settitle{\begin{center}\baselineskip14\p@\relax\scshape\large\@title\end{center}}
\def\@oddhead{\baselineskip14\p@\relax\thepage\hfil\scshape\headtitle\hfil}
\def\@evenhead{\baselineskip14\p@\relax\hfil\scshape\sectitle\hfil\thepage}
\makeatother

\title{\headtitle}
\author{\mauthor}
\date{\today}

\begin{document}

\sloppy
\allowdisplaybreaks

\maketitle



For as long as I can remember, almost every foreword to any decent math olympiad book begins with mentioning that the very first official math contest for high school students was held in Hungary in~$1894$. Namely, that was the year when the graduates of all high schools and gymnasiums in Hungary were offered a chance to take part in a math contest organized by Hungarian Mathematical and Physical Society together with the popular magazine \txtsc{``Közé\-pis\-kolai Mate\-ma\-ti\-kai és Fizi\-kai Lapok''}.\ftnnote{The name translates as the ``Journal of Mathematics and Physics for Secondary Schools''.} Over the entire country, the event was attended by $54$ graduates, with $29$ of them submitting their work. Among those, only eight correctly solved even one problem (two contestants solved all three).\ftnnote{This info is taken from article~\cite{aaw}. See also books \cite{HMO_I},~\cite{HMO_II}.} And so began the famous \infolink{Eötvös-Kürschák contest}{https://imomath.com/index.cgi?page=collectionHun}.\ftnnote{This name was actually never officially used. Until $1943$ this competition was called Eötvös contest in dedication to a distinguished Hungarian physicist \infolink{Lorand Eötvös}{https://en.wikipedia.org/wiki/Lorand_Eotvos}, who at the time chaired the Hungarian Mathematical and Physical Society; after the World War II it was revived as the Kürschák competition to honor a well-known Hungarian mathematician and educator \infolink{József Kürschák}{https://en.wikipedia.org/wiki/Jozsef_Kurschak}.}

It bears noticing that a somewhat similar contest was already being held (since $1887$) in the Russian Empire. In $1884$, the professor of Kiev University \infolink{\Fio(Vassily P.Erma\-kov)}{https://ru.wikipedia.org/wiki/\%D0\%95\%D1\%80\%D0\%BC\%D0\%B0\%D0\%BA\%D0\%BE\%D0\%B2,_\%D0\%92\%D0\%B0\%D1\%81\%D0\%B8\%D0\%BB\%D0\%B8\%D0\%B9_\%D0\%9F\%D0\%B5\%D1\%82\%D1\%80\%D0\%BE\%D0\%B2\%D0\%B8\%D1\%87} became the publisher and the chief editor of ``\txtsc{The Journal of Elementary Mathematics}''. In $1886$, this magazine was transformed into another one called \infolink{\txtsc{``The Herald of Experimental Physics and Elementary Mathematics}''}{https://ru.wikipedia.org/wiki/\%D0\%92\%D0\%B5\%D1\%81\%D1\%82\%D0\%BD\%D0\%B8\%D0\%BA_\%D0\%BE\%D0\%BF\%D1\%8B\%D1\%82\%D0\%BD\%D0\%BE\%D0\%B9_\%D1\%84\%D0\%B8\%D0\%B7\%D0\%B8\%D0\%BA\%D0\%B8_\%D0\%B8_\%D1\%8D\%D0\%BB\%D0\%B5\%D0\%BC\%D0\%B5\%D0\%BD\%D1\%82\%D0\%B0\%D1\%80\%D0\%BD\%D0\%BE\%D0\%B9_\%D0\%BC\%D0\%B0\%D1\%82\%D0\%B5\%D0\%BC\%D0\%B0\%D1\%82\%D0\%B8\%D0\%BA\%D0\%B8} with another Kievan enthusiast of extracurricular science and math, teacher, and journeyman mathematician, \infolink{\Fio(Emmanuil K.Shpa\-chin\-sky)}{https://ru.wikipedia.org/wiki/\%D0\%A8\%D0\%BF\%D0\%B0\%D1\%87\%D0\%B8\%D0\%BD\%D1\%81\%D0\%BA\%D0\%B8\%D0\%B9,_\%D0\%AD\%D1\%80\%D0\%B0\%D0\%B7\%D0\%BC_\%D0\%9A\%D0\%BE\%D1\%80\%D0\%BD\%D0\%B5\%D0\%BB\%D0\%B8\%D0\%B5\%D0\%B2\%D0\%B8\%D1\%87}). Every issue of the magazine (they were published twice a month!) had a prob\-lem-sol\-ving section that contained several problems to be attempted by the readers, with the annual contest winners rewarded by books on mathematics and physics. In $1887$, this new magazine also started including \infolink{problems for the high school students}{https://www.vofem.ru/ru/issues/1887/3/8/32/}.\ftnnote{A short review of the history of those two remarkable magazines can be found in~\cite{Busev}.}

Despite some similarities---both these contests targeted the same audience (the high school students and gymnasium graduates), and both were publicized and delivered via a popular scientific magazine---they also had several very significant differences.

\vmskip
\centerline{%
\pic{0.55}{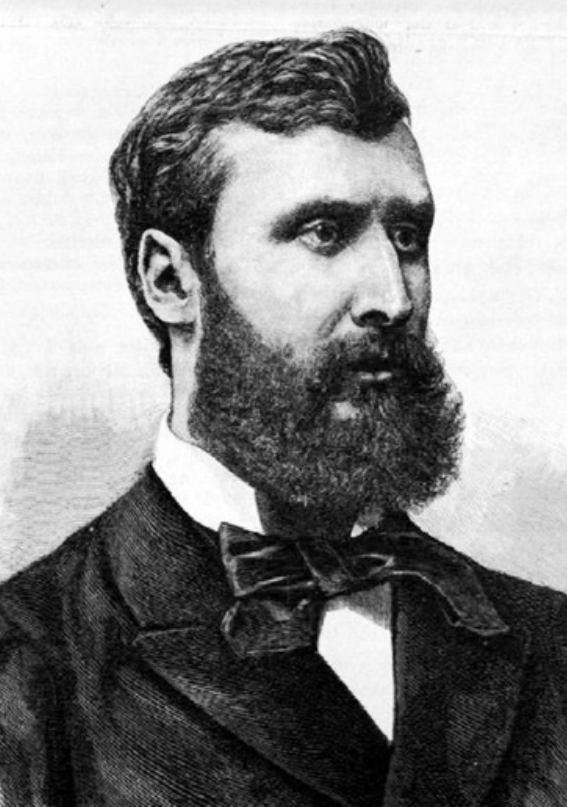}%
\qquad\qquad
\pic{0.75}{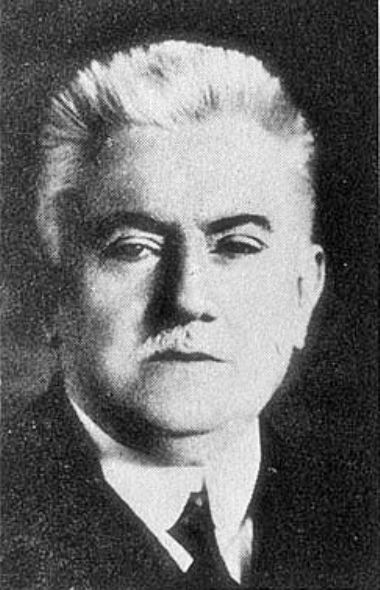}%
}
\ccaption{\hspace{4mm}Lorand Eötvös ($1848$--$1919$)\hspace{12mm} József Kürschák ($1864$--$1933$)}

First, the Eötvös contest had a formal jury, consisting of professional mathematicians and teachers. Their job was to create the problem list, but also they were the ones to read and evaluate the contest papers. Quite differently, the problems for the Herald contest were sent to the editor by the readers, and the evaluation was done by the editor himself with some help from a few of his equally enthusiastic collaborators.

Second, the Eötvös competition was held via a very precise and formal framework---all the problem solving had to be done on a specific day (under supervision by the math teachers from local gymnasiums) and then mailed to Budapest before a set date; there were always three problems and there had to be two prizes awarded (the winner and the runner-up). As for the Herald contest, the schoolkids solved the problems at home, doing that over at least several days or even weeks.

Third, these contests had a very meaningful difference in accessibility. Not every gymnasium graduate could take part in the Eötvös competition, due to some formal requirements for the local organization of the contest. The Herald contest was more accessible, but still a participant had to somehow get their hands (by buying or borrowing) on a current issue of the magazine.

Fourth, the Herald contest was quite informal, very similar to the \infolink{\txtsc{``Kvant''}}{https://en.wikipedia.org/wiki/Kvant_(magazine)} magazine problem solving contest.\ftnnote{That contest started in $1970$ when \txtsc{``Kvant''} magazine was founded.} There were no official awards; at the end of each year, the publisher tallied up the results and sent some popular science books to the best solvers. In contrast, the Eötvös contest had a very formal structure coupled with a substantial official component---the winner would receive a prize awarded by Hungarian Society of Science and Mathematics together with a significant sum of money and an invitation to enroll at the School of Science of \infolink{University of Budapest}{https://en.wikipedia.org/wiki/Eotvos_Lorand_University}. In the later years, when the number of participants has noticeably increased, this invitation (or at least, the recommendation) was awarded to the students who authored the ten best papers.

We should conclude that only one of these two competitions---obviously, we mean the Eötvös contest---qualifies to be called a math olympiad (even if it was done by correspondence).

\vmskip
\centerline{%
\pic{0.695}{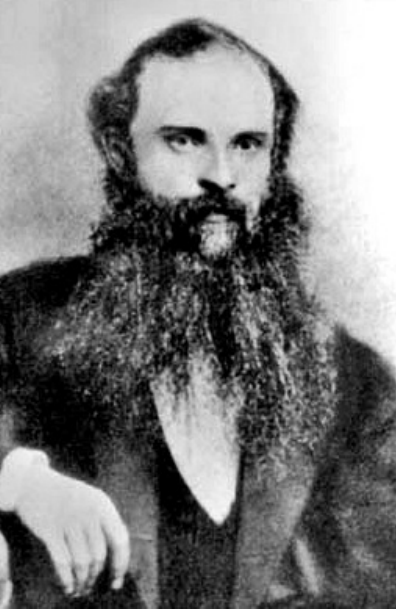}%
\qquad\qquad
\pic{1.751}{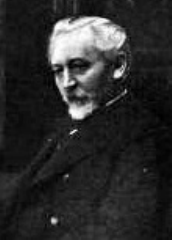}%
}
\ccaption{\hspace{2mm}\fior(V.P.Ermakov) ($1845$--$1922$)\hspace{10mm}\fior(E.K.Shpachinsky) ($1848$--$1912$)}

It bears emphasizing that the \txtsc{``Herald''}, though based in Kiev and then in Odessa, was by no means a provincial publication. It had a truly All-Russian span and coverage. Article \cite{Gushel} takes as a typical example issue \tnumsign{560} printed in~$1912$. Questions for that issue came from Odessa, Samara, and Stavropol, while the solutions were submitted by the readers from \stpb, Moscow, Kazan, Astrakhan, Warsaw, Odessa, Armavir, Łódź, and Sterlitamak. The popularity of the \txtsc{``Herald''} in the country was enough that even many years after the magazine ceased to exist, its issues could still be found in the public libraries of almost every large Russian city.

Some Russian sources (for instance, in the article entitled \infolink{History of the olympiad movement}{https://olimpiada.ru/article/687/}) refer to some youth competitions allegedly held at the very end of XIX century by the \infolink{Russian Astronomical Society}{https://en.wikipedia.org/wiki/Russian_Astronomical_Society}. Such articles often referred to book \cite{lutsky}, which, however, contains no mention of the said contests. I became interested in the aforementioned claim and started digging a little deeper into this fascinating matter, trying to find some definitive reference. I also sent inquiries to some acquaintances involved in Russian astronomy olympiads.

Finally, one of my friends found the author of the original reference, \fio(Mikhail G.Gavrilov), who at some point was very active in All-Russia Astronomical Olympiads. Gavrilov claimed that he cannot now remember the source of that sentence but that he most likely read about it in some ancient book found on the shelves of an antique bookshop.\ftnnote{The most likely source would be either an issue of the Notices of the Russian Astronomical Society (which has been published since $1892$) or one of the annual almanacs of the Russian Astronomical Society (published by the society since $1909$).} This leaves us with the peculiar task of finding that reference. Until then the claim about the existence of some XIX century astronomy competitions for high school students remains suspect.

Regardless, my investigations were not wasted. Examining the alleged competitions of Russian schoolkids in the XIX century, I have, almost by sheer luck, run into a much more interesting (I would even dare say, shocking!) fact. It turned out that the very first\ftnnote{Unless some painstakingly meticulous researcher digs out and presents us with the proof of, say, Napoleonic era France holding official math contests (except the college entrance exams, of course).} mathematical contest for high school students, held quite formally and officially at a major city (or country) level, occurred not in Hungary but... where do you think? of course, in \stpb! Moreover, in this case we have more than enough documents and publications corroborating this claim---all of them easily accessible online.

\xxx

\textit{...Once upon a time, in a kingdom far far away, its good sovereign decided to hold the great festivities---but not the usual kind. Declared the king that he was bored with the ceaseless harp stringing and dragon wranglers; no longer did he desire to see the knights competing in archery and fencing, or engaging in deathly fights on a muddy field. Instead, he ordered to promptly gather the youths most clever and bright so that those who exceed the others with their wit and wisdom can be found and properly rewarded...}

\vsskip

Joking aside, here is what happened. Almost two hundred years ago, in the year $1840$, the superintendent of \stpb Educational District, all of a sudden decreed to hold a ``\textit{competition event}'' between the students of the city's gymnasiums (using more contemporary terminology, the high school students).

And here is what was written on that matter by a chronicler monk... beg your pardon, we meant Nikolay Anichkov, the humble historian of the Third \stpb Gymnasium, whose chronicle was published in the imperial capital merely $150$ years ago, in $1873$.\ftnnote{Here and later in this chapter, we will translate the older terminology in our quotes and excerpts, using current names and terms.}

\begin{qte}

\tquote{... during the annual tests of $1841$, there was an order from the superintendent of the district, Prince \fior(G.P.Volkonsky), to organize a contest between the students of all gymnasiums. Accordingly, the questions were selected and offered: in Russian language for three students from Grade VII, in Latin for the same number of Grade VI students, in history for the Grade V, and in algebra for Grade IV. All the students gathered at the university, and there, under the supervision of the district inspector \fior(P.P.Maximovich), wrote their papers without any help or textbooks, excluding dictionaries for the translations. Thus, written papers or translations were then examined by the [University] professors, and by their determination, the best submissions were awarded by books and diplomas. Four prizes were assigned to each subject, and so the number of the annual awards was~$16$. Such contests were held for three years $(1841$--$1843)$, and therefore, $48$ prizes were awarded, out of which $19$ were given to the $3$rd Gymnasium. Altogether our Gymnasium has sent $40$ students, and consequently, half of them were honored with the prizes.}

{\hfill \fior(N.M.Anichkov), \cite{Ist3gym}, pp.$156$--$157$}

\end{qte}

These recollections, however, contain two major errors (that fact is confirmed by various other chronicles and documents). The first is that the aforementioned contests began in $1840$, that is, one year earlier than reported by \fior(N.M.Anichkov). The second mistake was made in the name of the then superintendent of the district (in the current terminology, that would be the secretary of education for the Northwest region of the Russian Empire).

That important position in the year $1840$ was occupied by Prince Mikhail Alexandrovich \infolink{Dondukov-Korsakov}{https://ru.wikipedia.org/wiki/\%D0\%94\%D0\%BE\%D0\%BD\%D0\%B4\%D1\%83\%D0\%BA\%D0\%BE\%D0\%B2-\%D0\%9A\%D0\%BE\%D1\%80\%D1\%81\%D0\%B0\%D0\%BA\%D0\%BE\%D0\%B2,_\%D0\%9C\%D0\%B8\%D1\%85\%D0\%B0\%D0\%B8\%D0\%BB_\%D0\%90\%D0\%BB\%D0\%B5\%D0\%BA\%D1\%81\%D0\%B0\%D0\%BD\%D0\%B4\%D1\%80\%D0\%BE\%D0\%B2\%D0\%B8\%D1\%87}, who nowadays is usually remembered because of a rather raunchy epigram written by \infolink{Alexander Pushkin)}{https://en.wikipedia.org/wiki/Alexander_Pushkin}. The great poet was very much upset by the censorial oversight (Prince Dondukov was also the imperial capital's main censor, supervising all the publications printed in the city), which, from his point of view, contradicted the promise given to him by Emperor \infolink{Nicholas~I}{https://en.wikipedia.org/wiki/Nicholas_I_of_Russia} of Russia.

In the later years, when Pushkin was personally introduced to the prince, he ``upgraded'' his opinion of him (incidentally, the prince's younger brother \infolink{Nikolay Korsakov}{https://ru.wikipedia.org/wiki/\%D0\%9A\%D0\%BE\%D1\%80\%D1\%81\%D0\%B0\%D0\%BA\%D0\%BE\%D0\%B2,_\%D0\%9D\%D0\%B8\%D0\%BA\%D0\%BE\%D0\%BB\%D0\%B0\%D0\%B9_\%D0\%90\%D0\%BB\%D0\%B5\%D0\%BA\%D1\%81\%D0\%B0\%D0\%BD\%D0\%B4\%D1\%80\%D0\%BE\%D0\%B2\%D0\%B8\%D1\%87}, who died very young, was one of Pushkin's classmates in Tsarskoselsky Lyceum).

\wpic{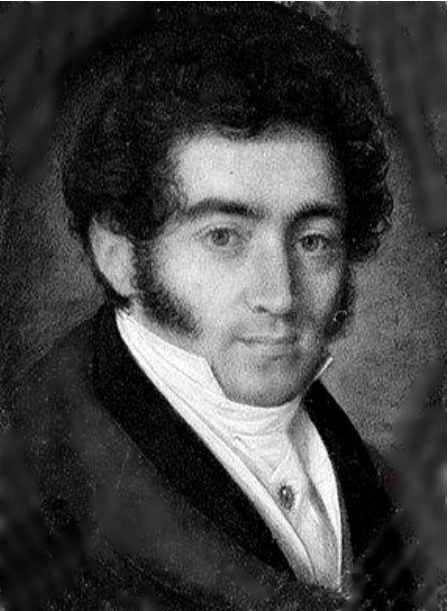}{s=0.72,w=0.35,c=\fior(M.A.Dondukov)}

Is it possible that the prince's most noble brain was struck by this ingenious idea? It is far more likely that the proposal came from the professors of the imperial university. Still, we should give to Prince Dondukov his due---he supported that ``revolutionary'' initiative, and, by the looks of it, took all necessary measures supporting this experiment. Generally speaking, Prince Mikhail Alexandrovich cannot be called a mindless government functionary---for instance, it was he who started the process of revision and unification of the gymnasial curricula of that era. By his order, starting in $1837$, the principals of all gymnasiums in the Northwest were directed to compile very detailed teaching plans on every subject taught at their schools and then to send them to the office of the superintendent. These documents were then examined, reviewed, and summarized by the Academic Council of the \stpb University. Over the next few years, the council accumulated all the suggestions and corrections made by the professors. Finally, in $1842$, the council produced the summary recommendations for each subject which were then ``taken into consideration'' when the gymnasial curricula were updated accordingly.

Here is another quote (this time, from the chronicles of the Second \stpb Gymnasium).

\begin{qte}

\tquote{During minister Uvarov's tenure,\ftnnote{\infolink{Sergey S.~Uvarov}{https://en.wikipedia.org/wiki/Sergey_Uvarov} ($1786$--$1855$) was an influential Russian imperial statesman, minister of education during the reign of Emperor Nicholas I, president of \stpb Academy of Sciences.} intending to \underline{excite} the competition between not only the students but also between teachers and schools, and to make the parents and the public aware of the successes of our youth, were held the contests for the students of Petersburg gymnasiums.

In $1840$, according to the recommendation of the Minister of Education from May $2$ of the same year, as an experiment, a common contest was established for the students of \stpb's four gymnasiums. The rules were set by the Minister as follows: the contest is held for the students of the four senior grades in common subjects taught to these years. Thus, the principal of each gymnasium selects, at his discretion, three students from each grade, who, at some fixed day and time, gather in one of the said gymnasiums. A common subject is chosen for each of these grades, e.g., Russian literature for Grade VII, and then the students are given three or four hours to write their papers, under supervision of a teacher, on a topic determined by the superintendent. The papers are then submitted to the superintendent, who, after replacing the name of a contestant with a corresponding number, transfers the papers to an unaffiliated and unbiased professor or another adjudicator to be evaluated according to their merits. After receiving the marks, the superintendent determines the winners in each subject; no more than four per grade, with two of them receiving books and two receiving commendation certificates.}

{\hfill \fior(A.V.Kurganovich), \fior(A.O.Kruglyi), \cite{Ist2gym}, pp.$371$--$372$}

\end{qte}

Similar notes mentioning these contests for the city's high school students appear in other chronicles of the capital's gymnasiums dated by the first half of the XIX century. They were also written up by \infolink{Georg Karl Schmid}{https://ru.wikipedia.org/wiki/\%D0\%A8\%D0\%BC\%D0\%B8\%D0\%B4,_\%D0\%93\%D0\%B5\%D0\%BE\%D1\%80\%D0\%B3_\%D0\%9A\%D0\%B0\%D1\%80\%D0\%BB} in his book surveying Russian secondary education of that epoch:

\begin{qte}

\tquote{Here we should note the rule prescribed in $1840$--$1842$ for \stpb Educational District stating that four best students of Grades IV--VII from all gymnasiums of the capital would gather at the university and sit, supervised by a district inspector, for an annual test in various subjects. Their papers were then submitted to the university professors for an evaluation, with the best students being awarded books and commendations.

...

The contest would be concluded with the public ceremony, where the general report on the state of all \stpb gymnasiums was read.}

{\hfill \fior(G.K.Schmid), \cite{Schmid}, pp.$320$--$321$}

\end{qte}

Now we can conclude that this is nothing else but an \textbf{official in-person}\ftnnote{In contrast to such cor\-res\-pon\-den\-ce-ba\-sed competitions as Eötvös-Kürschák or \txtsc{``Herald'' contests}.} \textbf{all-city olympiad} open to all qualified high school students of the Russian imperial capital, city of \stpb. Moreover, the rules and logistics of this contest are practically the same as the contemporary ones used for similar level competitions around the world.

Summarizing the available sources (see books \cite{Ist1gym}, \cite{Ist2gym}, \cite{Ist3gym}, and \cite{Schmid}), we can claim that the contests were held for three years straight, in $1840$--$1842$. At the very end of the school year, each gymnasium would send their team, consisting of three (or four, in $1842$) participants per each of the four subjects.

After the papers were evaluated, the four best results would receive the respective prizes. The winner and the runner-up were awarded some scientific and mathematical books, while the third and the fourth places were presented with a certificate of commendation.

In $1840$, the subjects were: \underline{algebra}, world history, Latin language, and Russian literature. The competition was held in the Second Gymnasium on June $4$ and $7$ (the first day---for Grades IV, V, and VI; the second day---for Grade VII). The awards were presented on June $22$ at the university (all dates here are given in the old Julian calendar, which was then still commonly used in Russia).

In $1841$, the contest was held in \underline{algebra}, world history, Latin language, and Russian history. It took place on the main campus of \stpb University on June $5$ and~$6$, and the award presentation was done there as well on June $23$.

In $1842$, the list of subjects consisted of \underline{geometry}, Latin and Greek languages, and Russian literature. The competition was once again held on the grounds of \stpb University on June $4$ and~$5$.

From our point of view, it is important to note that the four competition subjects always included a mathematical discipline. Namely, Grade IV students (approximately $14$--$15$ years old, which corresponds to Grade $9$ of contemporary school) wrote their papers, solving problems in algebra and geometry; Grade V students competed in world history, Grade VI in Latin, and Grade VII in Russian literature or history.

Altogether in three years $48$ prizes were distributed---three years, four subjects, and four prizes per subject.

Out of these $48$ diplomas, one was gained by the students from the First Gymnasium, $20$ by the students of the Second Gymnasium, $19$ by the students of the Third Gymnasium, and accordingly, $8$ by the students of the Fourth (Larin) Gymnasium.

Unfortunately, the names of the laureates are known only for the First and the Second Gymnasiums. For the Fourth Gymnasium, we know the names only for the very first contest.

\begin{qte}

\tquote{In $1840$, by the order of the Minister of Education, a common contest was held for the students of four senior classes of all gymnasiums. In that contest the students of Larin Gymnasium distinguished themselves, especially in the subjects of history and mathematics. Student Schedrin I$^{st}$ and Moshnin received the first and the second prize in history; student Grigoriev (later a teacher and an inspector at the same gymnasium) was awarded the third prize in mathematics.}

{\hfill \fior(K.M.Blumberg), \fior(V.P.Ostrogorsky), \cite{Ist4gym}, p.$18$}

\end{qte}

In the gymnasium chronicles we can read about the details on organization and conduct of these city-level contests---down to the little things like special paper to be used by the contestants (the stamped sheets were provided by the superintendent's office); likely, this was done to prevent possible cheating.

\begin{qte}

\tquote{Topics and questions for the competition were delivered sealed and were opened only after all selected students assembled together. In the Second Gymnasium, the contest was held on June $4$ and $7$; the students wrote their papers from $10\nsfrac12$~a.m.\ to $12\nsfrac12$~p.m.\ and then after lunch until $5$~p.m., whereupon the principal collected the submissions and sent them sealed to the superintendent. The following students took part:

In algebra, Grade $4$: Ivan Glazunov, Fyodor Mednikov, Konstantin Timofeev.

...

Upon completion of the contest, the formal act of the four gymnasiums was held on June $22$ in the presence of the top management and the public. Among the students of the Second Gymnasium, the following were honored with the awards: in algebra, the second prize went to Grade IV student Fyodor Mednikov, and the fourth prize---to Konstantin Timofeev, student of the same year.}

{\hfill \fior(A.V.Kurganovich), \fior(A.O.Kruglyi), \cite{Ist2gym}, p.$31$}

\end{qte}

It seems that the best results in mathematics were achieved by the Third Gymnasium. It is most unfortunate that the chronicle of that particular school mentions only the general statistics but does not provide any names.

\xxx

Now, we will dig a little deeper into the details. Let us find out what exactly was the \infolink{\stpb Educational District}{https://ru.wikipedia.org/wiki/\%D0\%A1\%D0\%B0\%D0\%BD\%D0\%BA\%D1\%82-\%D0\%9F\%D0\%B5\%D1\%82\%D0\%B5\%D1\%80\%D0\%B1\%D1\%83\%D1\%80\%D0\%B3\%D1\%81\%D0\%BA\%D0\%B8\%D0\%B9_\%D1\%83\%D1\%87\%D0\%B5\%D0\%B1\%D0\%BD\%D1\%8B\%D0\%B9_\%D0\%BE\%D0\%BA\%D1\%80\%D1\%83\%D0\%B3} in the $1840$s, how many gymnasiums the capital had, what was their span and make-up---namely, how many students they taught, from which social strata of the population they came etc.

The educational district centered around the capital comprised the entire Northwest of the Russian empire. By the imperial decree of January $12$, $1831$, that district covered the \stpb \infolink{governorate}{https://en.wikipedia.org/wiki/Governorate_(Russia)} as well as Novgorod, Pskov, Vologda, Olonets, and Arkhangelsk governorates. Until $1835$ the district was managed by the so-cal\-led \textit{college} committee at the \stpb University. However, beginning on June $25$, $1835$, when the new regulation for educational districts was adopted, virtually all governance responsibilities were transferred to the office of the district's superintendent. From this moment on, the university administration (first of all, its rector) was only in charge of education at the level of universities and colleges of the capital and its environs.

The superintendent of the \stpb Educational District at the time was the already mentioned by us Prince Mikhail Kor\-sa\-kov-Don\-du\-kov, executing those duties from April $30$, $1832$ until May $7$, $1842$. Therefore, the entire education in the capital was in the hands of two people---Prince Dondukov and the rector of the imperial university, \infolink{his excellency}{https://en.wikipedia.org/wiki/Active_State_Councillor} \infolink{Pyotr Pletnyov}{https://en.wikipedia.org/wiki/Pyotr_Pletnyov}, literary critic and poet, formerly a close friend of Alexander Pushkin, and professor of Russian literature, who went on to occupy the top university post for twenty years. As the rector of the \stpb University, he also sat on the District Advisory Council, one of whose main tasks was improvement of the middle and high school system.

\wpic{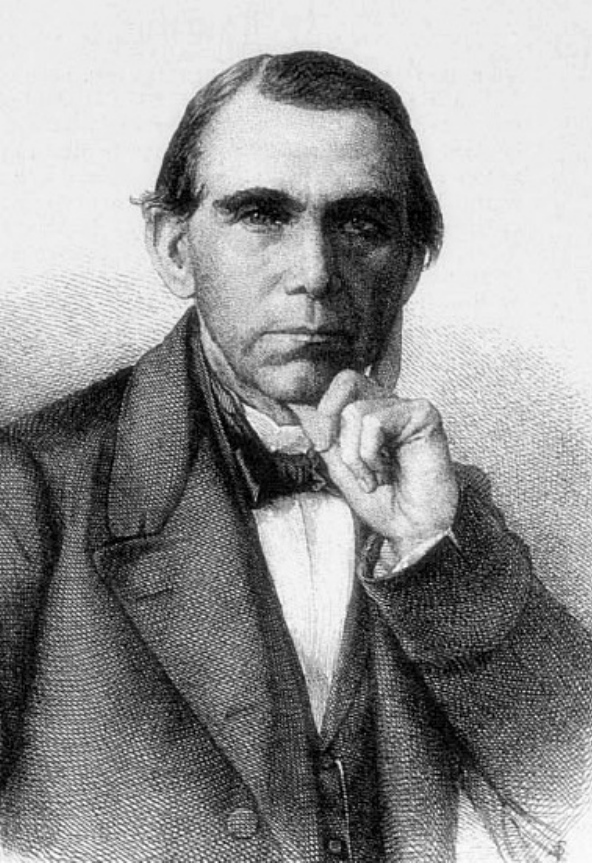}{s=0.55,w=0.35,c=\fior(P.A.Pletnyov)}

According to the contemporary memoirs, these two people respected each other, and their attitudes toward education were quite similar. Therefore, it seems plausible that the ``olympiad'' initiative came from Pletnyov himself or from some university professors around him. One of the chronicles actually mentions the rector's active participation in setting up the contest. It says that \glm{\textit{rector Pletnyov, who then managed the Petersburg educational district, asked Postels to keep an eye on this matter and to oversee the overall management of the contest}}.\ftnnote{See~\cite{Ist2gym}, p.$372$. An obvious mistake is made here; after $1835$ the rector no longer was in charge of the district. Mentioned above \infolink{\fior(A.F.Postels)}{https://en.wikipedia.org/wiki/Alexander_Postels} served then as the principal of the Second Gymnasium.}

For the sake of completeness, we should mention that in these years, the math at the \stpb university was taught by the chair of the pure and applied mathematics, full professor \infolink{\Fio(Dmitri S.Chizhov)}{https://ru.wikipedia.org/wiki/\%D0\%A7\%D0\%B8\%D0\%B6\%D0\%BE\%D0\%B2,_\%D0\%94\%D0\%BC\%D0\%B8\%D1\%82\%D1\%80\%D0\%B8\%D0\%B9_\%D0\%A1\%D0\%B5\%D0\%BC\%D1\%91\%D0\%BD\%D0\%BE\%D0\%B2\%D0\%B8\%D1\%87}, full professor \infolink{\Fio(Vikenty A.Ankudovich)}{https://ru.wikipedia.org/wiki/\%D0\%90\%D0\%BD\%D0\%BA\%D1\%83\%D0\%B4\%D0\%BE\%D0\%B2\%D0\%B8\%D1\%87,_\%D0\%92\%D0\%B8\%D0\%BA\%D0\%B5\%D0\%BD\%D1\%82\%D0\%B8\%D0\%B9_\%D0\%90\%D0\%BB\%D0\%B5\%D0\%BA\%D1\%81\%D0\%B0\%D0\%BD\%D0\%B4\%D1\%80\%D0\%BE\%D0\%B2\%D0\%B8\%D1\%87}, and associate professor \infolink{Fyodor V.Chizhov}{https://ru.wikipedia.org/wiki/\%D0\%A7\%D0\%B8\%D0\%B6\%D0\%BE\%D0\%B2,_\%D0\%A4\%D1\%91\%D0\%B4\%D0\%BE\%D1\%80_\%D0\%92\%D0\%B0\%D1\%81\%D0\%B8\%D0\%BB\%D1\%8C\%D0\%B5\%D0\%B2\%D0\%B8\%D1\%87}.\ftnnote{In fact, their formal titles were ``\infolink{professor ordinarius}{https://en.wikipedia.org/wiki/Academic_ranks_in_Germany\#Main_positions}'', ``\infolink{professor extraordinarius}{https://en.wikipedia.org/wiki/Academic_ranks_in_Germany\#Main_positions}'', and ``\infolink{adjunct professor}{https://en.wikipedia.org/wiki/Adjunct_professor}''; here we use their approximate contemporary analogs.} 


Interestingly enough, the undisputed leader of the Russian mathematics of that era, academician \infolink{Mikhail Ostrogradsky}{https://en.wikipedia.org/wiki/Mikhail_Ostrogradsky}, did not teach at the university---his schedule was quite full with lectures at the \infolink{Institute of the Corps of Railroad Engineers}{https://ru.wikipedia.org/wiki/\%D0\%98\%D0\%BD\%D1\%81\%D1\%82\%D0\%B8\%D1\%82\%D1\%83\%D1\%82_\%D0\%B8\%D0\%BD\%D0\%B6\%D0\%B5\%D0\%BD\%D0\%B5\%D1\%80\%D0\%BE\%D0\%B2_\%D0\%BF\%D1\%83\%D1\%82\%D0\%B5\%D0\%B9_\%D1\%81\%D0\%BE\%D0\%BE\%D0\%B1\%D1\%89\%D0\%B5\%D0\%BD\%D0\%B8\%D1\%8F}, \infolink{General Pedagogical Institute}{https://ru.wikipedia.org/wiki/\%D0\%93\%D0\%BB\%D0\%B0\%D0\%B2\%D0\%BD\%D1\%8B\%D0\%B9_\%D0\%BF\%D0\%B5\%D0\%B4\%D0\%B0\%D0\%B3\%D0\%BE\%D0\%B3\%D0\%B8\%D1\%87\%D0\%B5\%D1\%81\%D0\%BA\%D0\%B8\%D0\%B9_\%D0\%B8\%D0\%BD\%D1\%81\%D1\%82\%D0\%B8\%D1\%82\%D1\%83\%D1\%82}, and \infolink{General Engineering College}{https://ru.wikipedia.org/wiki/\%D0\%9D\%D0\%B8\%D0\%BA\%D0\%BE\%D0\%BB\%D0\%B0\%D0\%B5\%D0\%B2\%D1\%81\%D0\%BA\%D0\%BE\%D0\%B5_\%D0\%B8\%D0\%BD\%D0\%B6\%D0\%B5\%D0\%BD\%D0\%B5\%D1\%80\%D0\%BD\%D0\%BE\%D0\%B5_\%D1\%83\%D1\%87\%D0\%B8\%D0\%BB\%D0\%B8\%D1\%89\%D0\%B5}.

Now let us return to the system of Russian gymnasial education of that time. In $1831$, gymnasiums in \stpb were all reformed into seven-year schools (until then, the length of education there was merely four years). For the first three years, the kids went through a basic gymnasium, district or church elementary school, and only after that, the most deserving (and most suitable, of course, as the social standing was certainly taken into account) boys\ftnnote{This is not a typo; until $1858$ only boys were allowed to attend the gymnasiums.} could transfer to the classical or technical gymnasium in order to obtain the more advanced school education.\ftnnote{The difference between the two types of gymnasiums was quite simple: the former included the so-called classical ``dead'' languages, Latin and Ancient Greek.}

The age of the gymnasiasts in each grade slightly fluctuated---sometimes the students were held back a year---but it is sufficient for us to know that the average age of the graduates (that is, the gymnasiasts who finished Grade VII and passed the exit exams) at the time was equal to approximately $17$--$18$ years. For instance, those who took part in the algebra contest (i.e., Grade IV students of \stpb gymnasiums) were about $14$ years old---that corresponds to the average age of a contemporary Grade $9$ student.

How many gymnasiums were there in the capital city of \stpb? As a matter of fact, they could be counted on one hand---in those years ($1840$--$1842$), the entire city was served by only four such schools (the Fifth Gymnasium was open in~$1845$).

Despite the population of the capital being around $470$ thousand, all four gymnasiums had just about nine hundred ($900$\,!) students.\ftnnote{The data is taken from \cite{IstStatSpb2},~p.$119$.} Incidentally, the entire Russian Empire at the time was inhabited by approximately $52$ million people.

The \infolink{First Gymnasium}{https://ru.wikipedia.org/wiki/\%D0\%9F\%D0\%B5\%D1\%80\%D0\%B2\%D0\%B0\%D1\%8F_\%D0\%A1\%D0\%B0\%D0\%BD\%D0\%BA\%D1\%82-\%D0\%9F\%D0\%B5\%D1\%82\%D0\%B5\%D1\%80\%D0\%B1\%D1\%83\%D1\%80\%D0\%B3\%D1\%81\%D0\%BA\%D0\%B0\%D1\%8F_\%D0\%B3\%D0\%B8\%D0\%BC\%D0\%BD\%D0\%B0\%D0\%B7\%D0\%B8\%D1\%8F} (formerly Noble Pansion) taught only the hereditary nobility; the \infolink{Second}{https://ru.wikipedia.org/wiki/\%D0\%92\%D1\%82\%D0\%BE\%D1\%80\%D0\%B0\%D1\%8F_\%D0\%A1\%D0\%B0\%D0\%BD\%D0\%BA\%D1\%82-\%D0\%9F\%D0\%B5\%D1\%82\%D0\%B5\%D1\%80\%D0\%B1\%D1\%83\%D1\%80\%D0\%B3\%D1\%81\%D0\%BA\%D0\%B0\%D1\%8F_\%D0\%B3\%D0\%B8\%D0\%BC\%D0\%BD\%D0\%B0\%D0\%B7\%D0\%B8\%D1\%8F} (formerly known as the Governorate Gymnasium) taught the sons of the government officials and minor nobility, who did not belong to the imperial court, as well as some wealthy merchants and various \infolink{service class people}{https://en.wikipedia.org/wiki/Raznochintsy}. The \infolink{Third}{https://ru.wikipedia.org/wiki/\%D0\%A2\%D1\%80\%D0\%B5\%D1\%82\%D1\%8C\%D1\%8F_\%D0\%A1\%D0\%B0\%D0\%BD\%D0\%BA\%D1\%82-\%D0\%9F\%D0\%B5\%D1\%82\%D0\%B5\%D1\%80\%D0\%B1\%D1\%83\%D1\%80\%D0\%B3\%D1\%81\%D0\%BA\%D0\%B0\%D1\%8F_\%D0\%B3\%D0\%B8\%D0\%BC\%D0\%BD\%D0\%B0\%D0\%B7\%D0\%B8\%D1\%8F} and the \infolink{Fourth}{https://ru.wikipedia.org/wiki/\%D0\%9B\%D0\%B0\%D1\%80\%D0\%B8\%D0\%BD\%D1\%81\%D0\%BA\%D0\%B0\%D1\%8F_\%D0\%B3\%D0\%B8\%D0\%BC\%D0\%BD\%D0\%B0\%D0\%B7\%D0\%B8\%D1\%8F} (aka Larin Gymnasium), theoretically speaking, accepted children of all social groups---but at the same time, the education was not free. As such, a son of a peasant, a factory worker, or a laundrywoman had a very low chance of enrolling there.

Moreover, the number of graduates was considerably lower than the total number of students. For example, in $1840$--$42$ all \stpb gymnasiums graduated about $40$ students per year! Just think how inadequate this number looks for the largest and the most influential city of a huge empire. The entire Northwest of Russia graduated about $80$--$100$ gymnasiasts per year, which is not really surprising since the entire enormous territory of \stpb, Novgorod, Pskov, Vologda, Olonets, and Arkhangelsk governorates (with a total population of $3.5$ million and total area of $1.25$ million square kilometers\ftnnote{Just for comparison: this is about the same as the total area of France, Germany, and Italy put together.}) was at the time served by only $9$\,(nine!) gymnasiums. And of course, the majority of those graduates came from the families of nobility and government officials.

\ccpic{0.35}{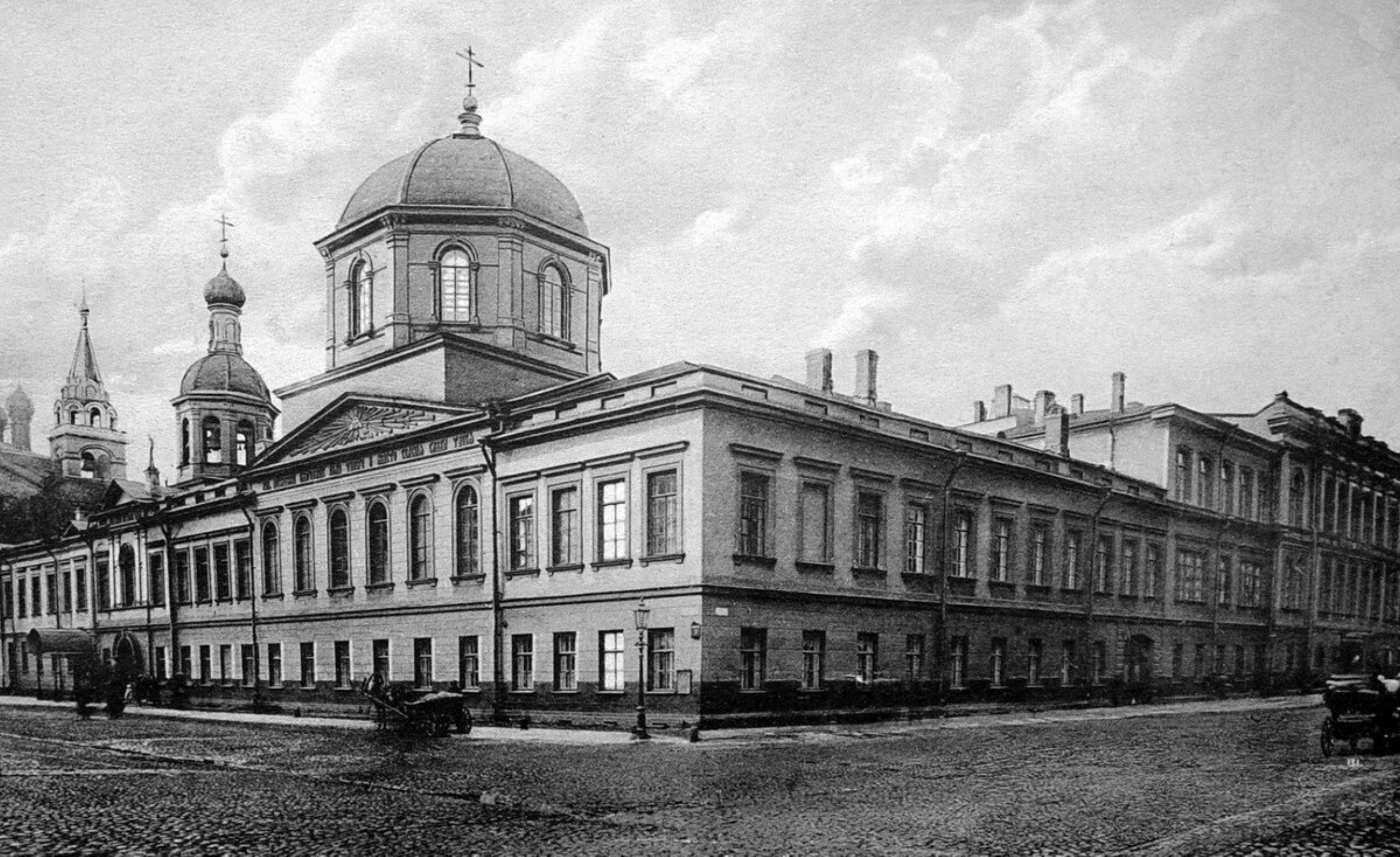}{The First \stpb Gymnasium ($7$ Ivanovskaya Street)}

We will continue with another interesting quote, this time from the chronicle of the First Gymnasium:

\begin{qte}

\tquote{... two students, Philip Depp and Nikolai Kaufman, were so successful in their studies on every subject, not excluding Greek language, that after the final tests, by determination of the Gymnasium council, they were honored by awarding the former student with the gold and the latter---with the silver medal. Such an event... was especially remarkable for the First Gymnasium because it was \textit{the very first time, since the Gymnasium was founded, these awards were conferred}. Furthermore, the first one of these laureates, in addition to the medal, was given another, even higher, recognition: in that year ($1840$), on the day of the final graduation act, by the order of the district authorities, student Depp received the second prize in Russian history at the award ceremony for the competition of all four \stpb gymnasiums on the subjects given in accordance with the proposal of the district superintendent, Prince Don\-du\-kov-Kor\-sa\-kov.}

{\hfill \fior(D.N.Soloviev), \cite{Ist1gym}, p.$195$}

\end{qte}

\ccpic{0.691}{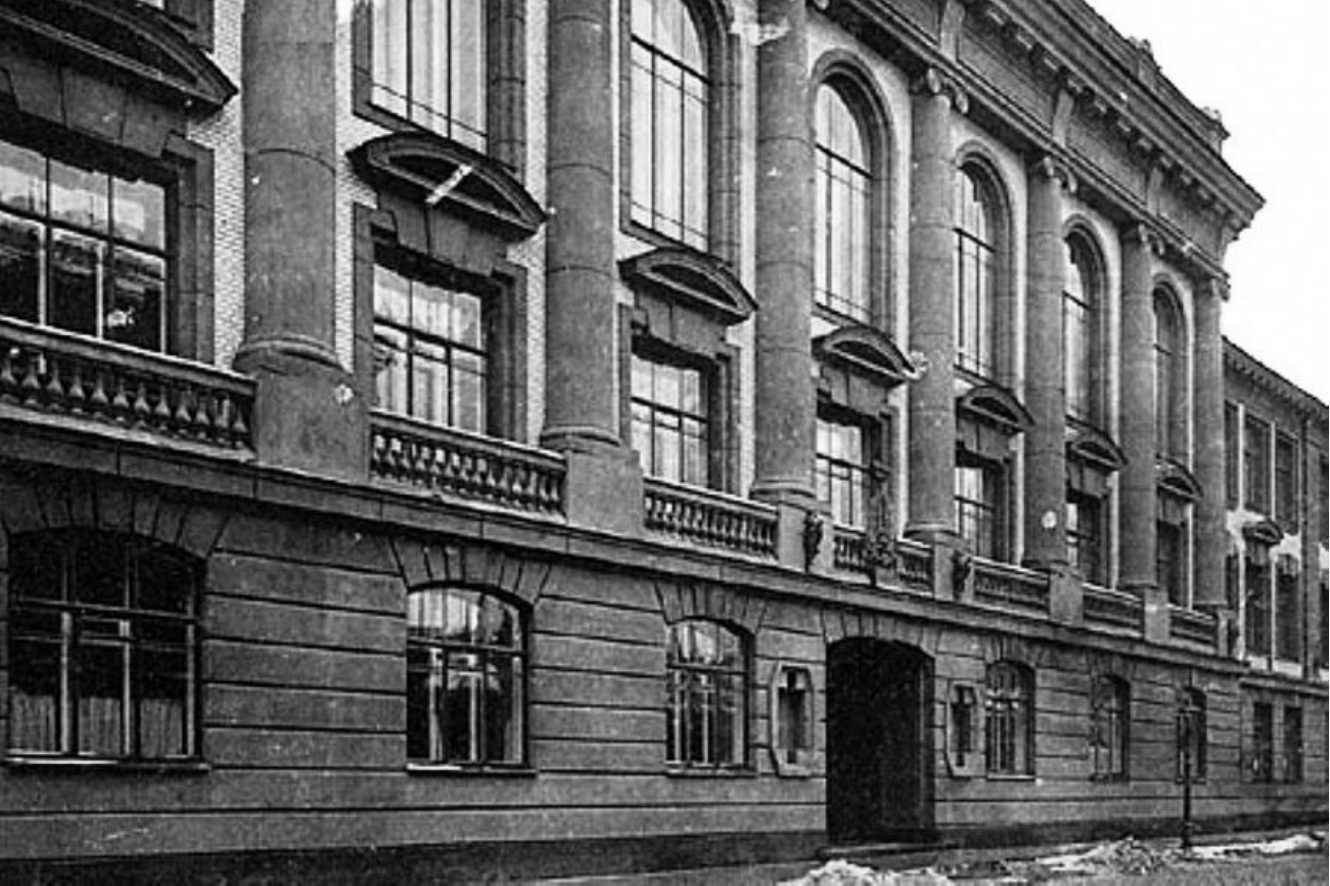}{The Second \stpb Gymnasium ($27$ Kazanskaya Street)}

\vmskip

\ccpic{0.474}{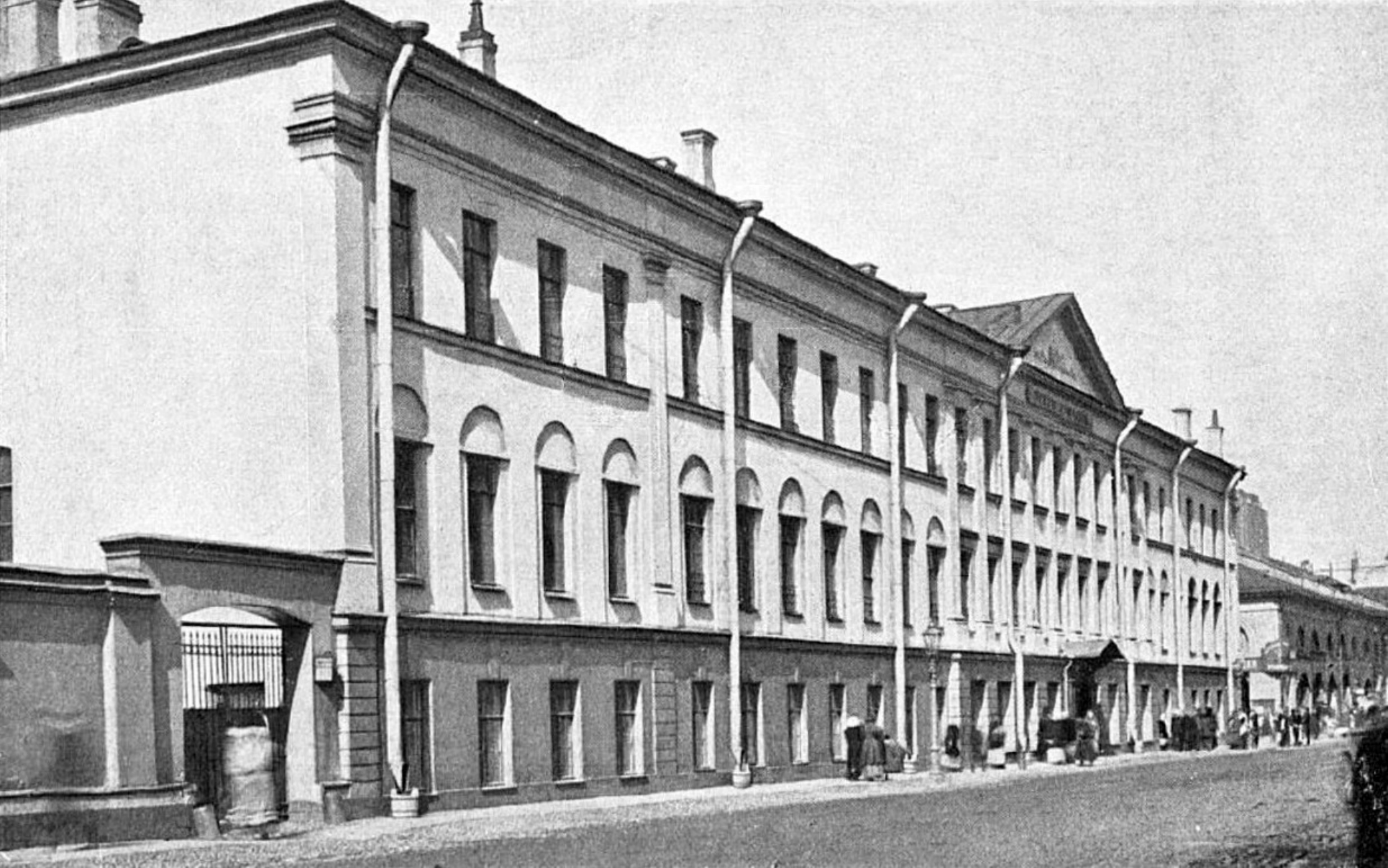}{The Third \stpb Gymnasium ($23$ Gagarinskaya Street)}

What exactly are they talking about here? First, it says that after the final graduation exams in the First Gymnasium in $1840$, two students received, respectively, a gold and a silver medal. Also, it says that in all previous years, none of the graduates of that school were awarded these decorations; therefore, it was quite a remarkable event for the gymnasium. But on top of all that, it turns out that the second place in the all-city olympiad in Russian history, taken by one of the students (Philip Depp), was deemed to be an even more prestigious award than his gold medal.

\ccpic{0.753}{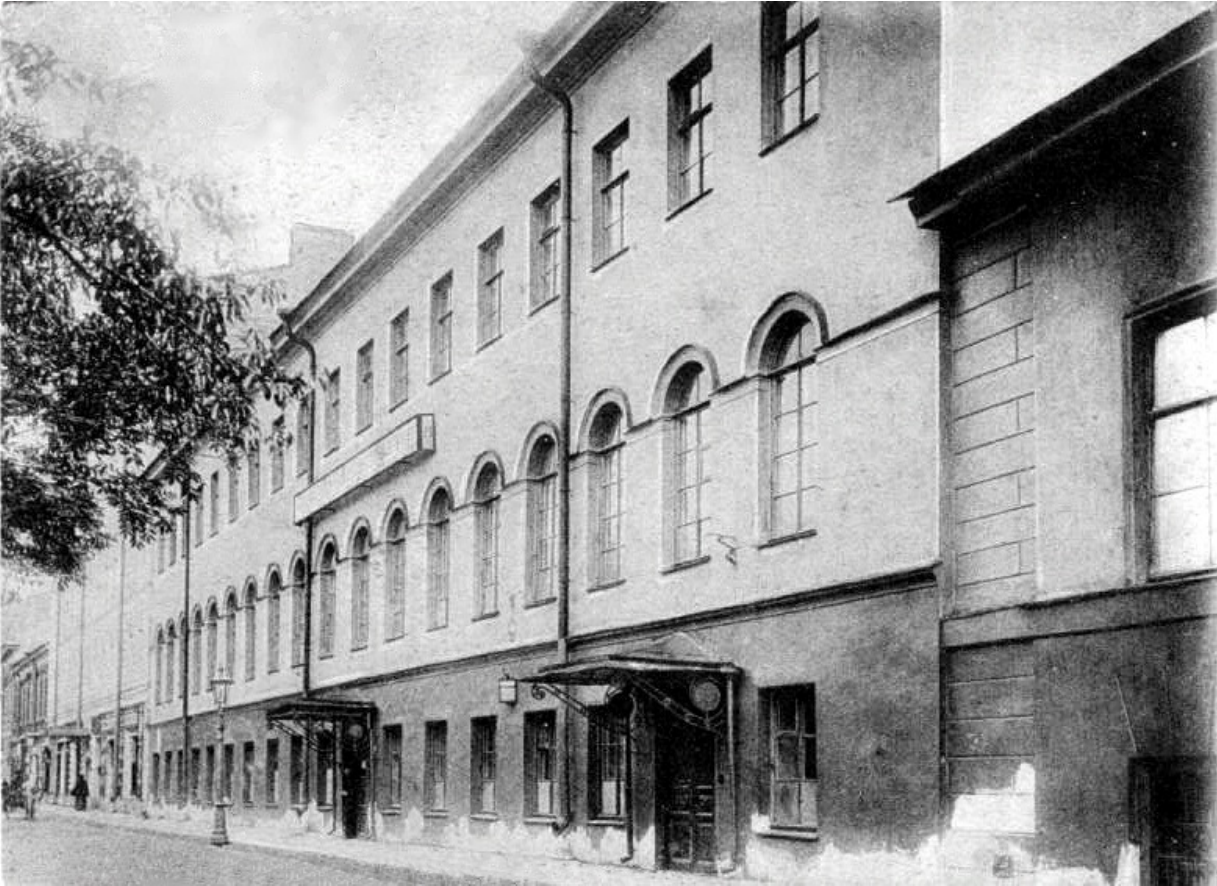}{The Fourth (Larin) \stpb Gymnasium ($15$, Sixth Line of Vassilievsky Island)}

Unfortunately, we know next to nothing about the future fate of all these contestants. That is quite regrettable---it would be, without a doubt, very interesting to learn what happened, say, to Pyotr Moszhechkov, a student of the Second Gymnasium and the triple laureate of these first olympiads. This ``wonder kid'' received the first prize in algebra in $1841$, while in $1842$ he managed to get the third prize in geometry and the first prize in Russian literature at the same time. 

One other multiple prize winner was \infolink{Konstantin Timofeev}{https://ru.wikipedia.org/wiki/\%D0\%A2\%D0\%B8\%D0\%BC\%D0\%BE\%D1\%84\%D0\%B5\%D0\%B5\%D0\%B2,_\%D0\%9A\%D0\%BE\%D0\%BD\%D1\%81\%D1\%82\%D0\%B0\%D0\%BD\%D1\%82\%D0\%B8\%D0\%BD_\%D0\%90\%D0\%BA\%D0\%B8\%D0\%BC\%D0\%BE\%D0\%B2\%D0\%B8\%D1\%87} ($1827$--$1881$), the $1843$ graduate of the same Second Gymnasium, who won an award at all three contests---the fourth prize in algebra in $1840$, the third prize in world history in $1841$, and, finally, the third prize in Russian literature in $1842$. After graduation he entered the School of History and Philology of \stpb University, from which four years later he received a diploma of candidate with a gold medal. Later he taught Russian literature and language at \infolink{\stpb Smolny Institute}{https://en.wikipedia.org/wiki/Smolny_Institute_of_Noble_Maidens} (as well as at two other colleges). He reached the position of Active State Councillor, published numerous articles and books in literary critique, pedagogy, and history of literature.

\begin{qte}

\tquote{Another true friend of ours was the class inspector Timofeev. He taught us Russian literature in Grade I. He loved us and the Russian poets and writers with all of his heart. He would recite Derzhavin, Zhukovsky, and Pushkin with tears in his eyes. To not learn his lessons or to dislike his favorite writers was to offer him the deepest insult, which was, of course, something that none of us would ever dare. To please him and to receive his praise, we would learn by heart entire long poems, with some of our students being quite good at reciting them aloud.}

{\hfill \fior(A.S.Eshevskaya), a student at Smolny Institute, \cite{sidv}, p.$207$}

\end{qte}

\BeginPic
\wpic{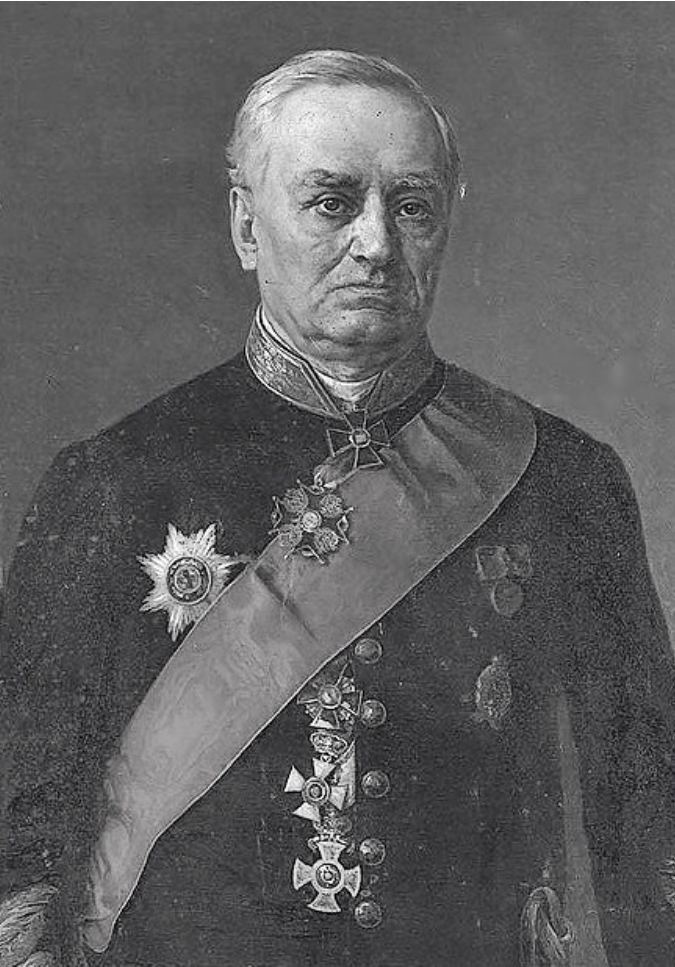}{s=0.44,w=0.35,c=\fior(I.I.Glazunov)}

Our final bright example---the fourth grader Ivan Glazunov, the member of the algebra ``team'' of the Second Gymnasium at the $1840$ contest. This was nobody else but the future mayor of \stpb, \infolink{Ivan Ilyich Glazunov}{https://ru.wikipedia.org/wiki/\%D0\%93\%D0\%BB\%D0\%B0\%D0\%B7\%D1\%83\%D0\%BD\%D0\%BE\%D0\%B2,_\%D0\%98\%D0\%B2\%D0\%B0\%D0\%BD_\%D0\%98\%D0\%BB\%D1\%8C\%D0\%B8\%D1\%87_(1826)}, a well-known Russian statesman (rising to the rank of the Privy Councillor), philanthropist, and businessman. The publishing company that he headed printed school textbooks, works of Russian educators, and classical Russian literature of the XIX century. He served as a mayor of the imperial capital in $1881$--$1885$; Ivan Ilyich also dedicated a lot of time and efforts to numerous social and charitable activities.

\EndPic

\xxx

So, let us summarize. More than $180$ years ago, \stpb's school district administration, together with the professors of the city university, held, for three years in a row, a city-level olympiad for high school students in four subjects (using the in-person written competition format). It seems clear that there was an intent to continue the contests indefinitely. Alas, the $1842$ contest turned out to be the last one---most likely, due to the change in the person of the superintendent. The new official, the former deputy of Prince \fior(M.A.Don\-du\-kov-Kor\-sa\-kov), Prince \infolink{Grigory Volkonsky}{https://ru.wikipedia.org/wiki/\%D0\%92\%D0\%BE\%D0\%BB\%D0\%BA\%D0\%BE\%D0\%BD\%D1\%81\%D0\%BA\%D0\%B8\%D0\%B9,_\%D0\%93\%D1\%80\%D0\%B8\%D0\%B3\%D0\%BE\%D1\%80\%D0\%B8\%D0\%B9_\%D0\%9F\%D0\%B5\%D1\%82\%D1\%80\%D0\%BE\%D0\%B2\%D0\%B8\%D1\%87}, obviously, did not demonstrate the same enthusiasm towards the idea of ``exciting'' of competitive spirit, exhibited by his predecessor. Without state support and administrative approval, the remarkable experiment withered away.

It is quite possible that the results of these first city ``olympiads'', as well as their questions and lists of participants, can still be found on the dusty shelves of the Central State Archive in \stpb. Right now, the author is unable to conduct the investigation in person---therefore, it remains my fervent hope that sooner or later someone will undertake and successfully finish the requisite historical and archival research.

\xxx

I would like to \underline{emphasize} that it is not really that important whether this or that specific country (or city) is formally recognized as the cradle of the first math olympiads. Still, to know your nation's history is not only interesting and useful, it is also necessary---at a minimum we should always strive to learn from the achievements (as well as the mistakes) of those who came before us.

\bigskip

\rightline{\txtsc{\Fio(Dmitri V.Fomin)}}

\vfill
\newpage


\def\bibauthor#1{\def\bAuthor{#1}}
\def\bibyear#1{\def\bYear{#1}}
\def\bibtitle#1{\def\bTitle{#1}}
\def\bibmisc#1{\def\bMisc{#1}}

\bibauthor{}\bibyear{}\bibtitle{}\bibmisc{}

\newcommand\bibentry[1]%
{
\bibitem{#1}%
\ifx\bAuthor\empty\else{\bAuthor,\ }\fi
{\textit{\bTitle}}%
\ifx\bMisc\empty,\ \else{.\ \bMisc,\ }\fi
{\bYear}\pp
}


\setcounter{section}{1}

%

\end{document}

\typeout{get arXiv to do 4 passes: Label(s) may have changed. Rerun}